\documentclass[11pt,a4paper]{article}
\usepackage[margin=2cm]{geometry}
\usepackage[english]{babel}
\usepackage[latin1]{inputenc}
\usepackage{amsmath, amsthm, amssymb}
\usepackage[pdftex]{graphicx}
\usepackage{amssymb}
\usepackage{latexsym}
\usepackage{amstext}
\usepackage{epstopdf}
\usepackage{floatflt}
\usepackage{hyperref}
\usepackage{enumitem}
\usepackage{multirow}
\usepackage{pdflscape}
\usepackage{makecell}
\usepackage{rotating}

\usepackage{eurosym}

\usepackage{slashed}

\usepackage{epstopdf}

\usepackage{nopageno}

\newtheorem{theorem}{Theorem}[section]

\newtheorem{definition}[theorem]{Definition}
\newtheorem{lemma}[theorem]{Lemma}
\newtheorem{proposition}[theorem]{Proposition}
\newtheorem{remark}[theorem]{Remark}

\newcommand{\aaa}{\mathfrak{a}}

\newcommand{\adc}{\aaa_2^c}

\newcommand{\adj}{\text{ad}}

\newcommand{\as}{\alpha}

\newcommand{\au}{\mathbf {a_1}}

\newcommand{\ba}{\begin{array}}

\newcommand{\bb}{\beta}

\newcommand{\be}{\begin{equation}}
\newcommand{\bea}{\begin{equation}\begin{array}}
\newcommand{\beal}{\begin{aligned}}
\newcommand{\beas}{\begin{equation*}\begin{array}}
\newcommand{\bef}{\begin{flalign}}
\newcommand{\befs}{\begin{flalign*}}
\newcommand{\bes}{\begin{equation*}}

\newcommand{\bit}{\begin{itemize}}
\newcommand{\blms}{{\mathfrak B}_{LS}}

\newcommand{\cc}{\mathbf{C}}

\newcommand{\dd}{\mathfrak D}
\newcommand{\ddn}{\mathbf {d_N}}
\newcommand{\ddnd}{\mathbf {d_{N-2}}}
\newcommand{\ddnt}{\mathbf {d_{N-3}}}

\newcommand{\dq}{\mathbf{d_4}}
\newcommand{\dqn}{\mathbf{d_{4n}}}

\newcommand{\ea}{\end{array}}

\newcommand{\eal}{\end{aligned}}

\newcommand{\ee}{\end{equation}}
\newcommand{\eea}{\end{array}\end{equation}}
\newcommand{\eeas}{\end{array}\end{equation*}}
\newcommand{\eef}{\end{flalign}}
\newcommand{\eefs}{\end{flalign*}}
\newcommand{\ees}{\end{equation*}}

\newcommand{\eit}{\end{itemize}}

\newcommand{\eo}{\mathbf{e_8}}

\newcommand{\eon}{\mathbf{e_8^{(n)}}}

\newcommand{\ep}{\varepsilon}

\newcommand{\es}{\mathbf{e_6}}
\newcommand{\esn}{\mathbf{e_6^{(n)}}}
\newcommand{\esquare}{\hfill $\square$}
\newcommand{\est}{\mathbf{e_7}}
\newcommand{\estn}{\mathbf{e_7^{(n)}}}

\newcommand{\fff}{\mathfrak F}

\newcommand{\fq}{\mathbf{f_4}}
\newcommand{\fqn}{\mathbf{f_4^{(n)}}}

\newcommand{\ggo}{\mathbf{\mathfrak g_0}}
\newcommand{\ggt}{\mathbf{\mathfrak g_{III}}}
\newcommand{\gh}{\gamma}

\newcommand{\gmf}{{\mathfrak g}}

\newcommand{\hri}{h_{\rho_i}}
\newcommand{\hrd}{h_{\rho_2}}
\newcommand{\hrt}{h_{\rho_3}}
\newcommand{\hru}{h_{\rho_1}}

\newcommand{\imp}{\ \Rightarrow\ }

\newcommand{\jotqq}{\mathbf{J_3^q}}

\newcommand{\joto}{\mathbf{J_3^8}}

\newcommand{\jp}{\circ}

\newcommand{\Ll}{\mathbb L}

\newcommand{\lms}{{\bf{\mathcal L_{MS}}}}

\newcommand{\lra}{\leftrightarrow}

\newcommand{\omv}{\omega^v}
\newcommand{\oms}{\omega^s}

\renewcommand{\proof}{\noindent {\bf Proof. }}

\newcommand{\pt}{\mathcal{P}}

\newcommand{\sref}[1]{{\bf\ref{#1}}}

\newcommand{\tfo}{T_O}
\newcommand{\tfop}{T_O^\prime}
\newcommand{\tfs}{T_S}
\newcommand{\tfsm}{T_S^-}
\newcommand{\tfsp}{T_S^+}
\newcommand{\tfspm}{T_S^\pm}
\newcommand{\tm}{T^-}
\newcommand{\tp}{T^+}

\newcommand{\um}{{\scriptstyle \frac12}}

\newcommand{\xd}{x_{P2}}
\newcommand{\xdb}{\bar x_{P2}}
\newcommand{\xpi}{x_{Pi}}
\newcommand{\xpj}{x_{Pj}}
\newcommand{\xpjb}{\bar x_{Pj}}

\newcommand{\xt}{x_{P3}}
\newcommand{\xtb}{\bar x_{P3}}
\newcommand{\xu}{x_{P1}}
\newcommand{\xub}{\bar x_{P1}}

\newcommand{\zz}{\mathbb Z}

\numberwithin{equation}{section}
\begin{document}

\begin{titlepage}
\begin{center}

\vskip 3.5cm

{\bf \huge Exceptional Periodicity\\\vskip 10pt and Magic Star Algebras\\\vskip 15pt \Large II : Gradings and {\em HT}-Algebras}

\vskip 1.5cm

{\bf \Large Piero Truini${}^{1,2}$, Alessio Marrani${}^{3,4}$, and Michael Rios${}^{5}$ }

\vskip 45pt

{\it ${}^1$ Quantum Gravity Research,\\ 101 S. Topanga Canyon Rd. 1159 Los Angeles, CA 90290 - USA}\\ \vskip 5pt

\vskip 25pt

{\it ${}^2$INFN, sezione di Genova,\\
Via Dodecaneso 33, I-16146 Genova, Italy}\\ \vskip 5pt

\vskip 25pt

 {\it ${}^3$Museo Storico della Fisica e Centro Studi e Ricerche ``Enrico Fermi'',\\
Via Panisperna 89A, I-00184, Roma, Italy}\\\vskip 5pt

\vskip 25pt

 {\it ${}^4$ Dipartimento di Fisica e Astronomia Galileo Galilei, Universit\`a di Padova,\\and INFN, sezione
di Padova, Via Marzolo 8, I-35131 Padova, Italy}\\\vskip 5pt

\vskip 25pt

{\it ${}^5$Dyonica ICMQG,\\5151 State University Drive, Los Angeles, CA 90032, USA}\\

\vskip 20pt

\texttt{truini@ge.infn.it},
\texttt{alessio.marrani@pd.infn.it},
\texttt{mrios@dyonicatech.com}

\end{center}

\vskip 75pt

\begin{center} {\bf ABSTRACT}\\[3ex]\end{center}


We continue the study of Exceptional Periodicity and Magic Star algebras,
which provide non-Lie, countably infinite chains of finite dimensional
generalizations of exceptional Lie algebras. We analyze the graded algebraic
structures arising in the Magic Star projection, as well as the Hermitian part of rank-3 Vinberg's  matrix algebras (which we dub {\em HT}-algebras), occurring on each vertex of the Magic Star.

\vskip 45pt







%
\vfill

\end{titlepage}

\newpage \setcounter{page}{1} \numberwithin{equation}{section}

\tableofcontents


\section{Introduction}

The aim of the present paper, the second of a series, is to continue the
rigorous formulation of Exceptional Periodicity. Such a formulation was
started in \cite{EP1}, in which we proved the existence of countably infinite, periodic chains of finite dimensional generalisations of the exceptional Lie algebras; in
particular, $\mathbf{e}_{8}$ has been shown to be part of a countably infinite family of algebras (named Magic Star algebras). By confining
ourselves to the simply laced cases, in \cite{EP1} we also determined the inner derivations and automorphisms of the Magic Star algebras.

The present paper focuses on the grading structures of the Magic Star algebras, in particular the 3-grading of $\estn$ and the 5-grading of $\eon$, and on the {\em HT}-algebras, which occur in each vertex of the Magic Star projection of Magic Star algebras, in the very same way cubic simple Jordan algebras occur in each vertex of the Magic Star projection of exceptional Lie algebras \cite{Truini, Mukai}.
Remarkably, {\em HT}-algebras were introduced some time ago by Vinberg in the study of homogeneous convex cones \cite{Vinberg}.
Endowed with a suitable commutative product, these algebras generalize cubic simple Jordan algebras, to which they reduce at the first $n=1$ level of Exceptional Periodicity, \cite{Truini, Mukai}.

\bigskip

The paper is organized as follows.

In Sec. \sref{s:lms} we recall the definition of the algebra $\lms$ introduced in \cite{EP1}. We summon up in particular the definition and main properties of the asymmetry function, lying at the core of the algebra product, that are extensively used in the present paper. In Sec. \sref{s:35} we analyze the 3-graded and 5-graded structures which can
be inferred by the shape of the Magic Star projection of Magic Star algebras.
Then, Sec. \sref{s:dta} deals with a generalization of cubic Jordan algebras given by {\em HT}-algebras, occurring on each of the six vertices of the Magic Star projection. Further structures named {\em HT}-pairs, which generalize the Jordan pairs \cite{loos1}, are introduced in Sec. \ref{s:pairs}. Finally, Sec. \ref{s:outlook} considers some future developments. An Appendix concludes the paper.


\section{The Magic Star algebra $\lms$}\label{s:lms}

We here recall the definition of the Magic Star algebra $\lms$, where $\lms$ is either $\esn$ or $\estn$ or $\eon$ (please see \cite{EP1} for the details).\\
We denote by $R$ the rank of $\lms$.
We recall that $N=4(n+1)$, $n=1,2,...$, hence $R=N-2=4n+2$ for $\esn$, $R=N-1=4n+3$ for $\estn$, $R=N=4n+4$ for $\eon$.\\
Let $V$ be a Euclidean space of dimension $R$ and $\{k_1 , ... , k_R\}$ an orthonormal basis in $V$. The set $\Phi$ of {\it generalized} roots of $\lms$ splits into $\Phi_O$ and $\Phi_S$, where:
\bea{l}
\Phi_O = \{ (\pm k_i \pm  k_j) \in \Phi \}
\\ \\
\Phi_S =  \{ \frac{1}{2} (\pm k_1 \pm k_2 \pm ... \pm k_N)  \in \Phi \}
\eea
We introduce the basis $\Delta = \{\alpha_1 , ... ,\alpha_R\}$ of $\Phi$, with $\alpha_i = k_i-k_{i+1}\, , \ 1\le i \le R-2$, $\alpha_{R-1} = k_{R-2}+k_{R-1}$ and $\alpha_R = -\frac12(k_1+k_2+...+k_N)$; we order $\Delta$ by setting $\alpha_i < \alpha_{i+1}$:
\be\label{sroots}
\Delta = \{k_1-k_2 < k_2-k_3< ...< k_{R-2} - k_{R-1}< k_{R-2} + k_{R-1}< -\um(k_1+k_2+...+k_N)\}
\ee

The Magic Star algebra $\lms$ is constructed as in \cite{EP1}:
\bit
\item[a)]  we select the set of ordered simple generalized roots $\Delta = \{\alpha_1 < ... <\alpha_R\}$ of $\Phi$
\item[b)] we select a basis $\{ h_1 ,...,h_R\}$ of the $R$-dimensional vector space $H$ over $\fff$ and set $h_\alpha = \sum_{i=1}^R c_i h_i$ for each $\alpha  \in \Phi$ such that $\alpha = \sum_{i=1}^R c_i \alpha_i$
\item[c)] we associate to each $\alpha  \in \Phi$ a one-dimensional vector space $L_\alpha$ over $\fff$ spanned by $x_\alpha$
\item[d)] we define $\lms = H \bigoplus_{\alpha \in \Phi} {L_\alpha}$ as a vector space over $\fff$
\item[e)] we give $\lms$ an algebraic structure by defining the following multiplication on the basis
$\blms = \{ h_1 ,...,h_R\} \cup \{x_\alpha \ | \ \alpha \in \Phi\}$, extended by linearity to a bilinear multiplication $\lms\times \lms\to \lms$:
	\be\begin{array}{ll}
	&[h_i,h_j] = 0 \ , \ 1\le i, j  \le R \\
	&[h_i , x_\alpha] = - [x_\alpha , h_i] = (\alpha, \alpha_i )\, x_\alpha \ , \ 1\le i \le R \ , \ \alpha \in \Phi \\
	&[x_\alpha, x_{-\alpha} ] = - h_\alpha\\
	&[x_\alpha,x_\beta] = 0 \ \text{for } \alpha, \beta \in \Phi \ \text{such that } \alpha + \beta \notin			 	\Phi \ \text{and } \alpha \ne - \beta\\
	&[x_\alpha,x_\beta] = \varepsilon (\alpha , \beta)\,  x_{\alpha+\beta}\ \text{for } \alpha , \beta \in \Phi \ \text{such that }  \alpha+ \beta \in \Phi\\
	\end{array} \label{comrel}\ee
\eit

where $\varepsilon (\alpha , \beta)$ is the {\it asymmetry function} defined as follows:
\begin{definition} Let $\Ll$ denote the lattice of all linear combinations of the simple generalized roots with integer coefficients
\be
\Ll = \left\{ \sum_{i=1}^R c_i \alpha_i \ |\ c_i \in \zz \ , \ \alpha_i \in \Delta \right\}
\label{lattice}
\ee
the asymmetry function $\varepsilon (\alpha , \beta): \ \Ll \times \Ll \to \{-1,1\}$ is defined by:
\be\label{epsdef}
\varepsilon (\alpha , \beta) = \prod_{i,j=1}^R \varepsilon (\alpha_i , \alpha_j)^{\ell_i m_j} \quad \text{for } \alpha = \sum_{i=1}^R \ell_i\alpha_i \ ,\ \beta = \sum_{j=1}^R m_j \alpha_j
\ee
where $\alpha_i , \alpha_j \in \Delta$ and
\be\label{epsdef1}
\varepsilon (\alpha_i , \alpha_j) = \left\{
\begin{array}{ll}
-1 & \text{if } i=j\\ \\
-1 & \text{if } \alpha_i + \alpha_j  \text{ is a root and } \alpha_i < \alpha_j\\ \\
+ 1 & \text{otherwise}
\end{array}
\right.
\ee
\end{definition}

The main properties of the asymmetry function are listed in the Propositions {\bf 5.1} and {\bf 5.2} of \cite{EP1}, whose statements we reproduce next:

\begin{proposition} \label{epsprop} The asymmetry function $\varepsilon$ satisfies, for $\alpha , \beta, \gamma , \delta \in \Ll$, $\alpha = \sum{m_i \alpha_i}$ and $\beta = \sum{n_i \alpha_i}$:
\bes\begin{array}{rrcl}
i) & \varepsilon (\alpha + \beta, \gamma) & =& \varepsilon (\alpha , \gamma)\varepsilon (\beta , \gamma) \\
ii) &\varepsilon (\alpha ,\gamma + \delta) & = &\varepsilon (\alpha , \gamma)\varepsilon (\alpha , \delta) \\
iii) &\varepsilon (\alpha , \alpha) & = &(-1)^{\frac12 (\alpha ,\alpha)-m_R^2 \frac{n-1}2}\\
iv) &\varepsilon (\alpha , \beta) \varepsilon (\beta , \alpha) &=& (-1)^{(\alpha , \beta) - m_R n_R (n-1)}\\
v) &\varepsilon (0 , \beta) &=& \varepsilon (\alpha , 0) = 1 \\
vi) &\varepsilon (-\alpha , \beta) &=& \varepsilon (\alpha , \beta)^{-1} = \varepsilon (\alpha , \beta) \\
vii) &\varepsilon (\alpha , -\beta) &=& \varepsilon (\alpha , \beta)^{-1} = \varepsilon (\alpha , \beta) \\
\end{array}\ees
\end{proposition}

\begin{proposition} If $\alpha, \beta ,\alpha+\beta \in \Phi$ then:
\beas{rll}
i) &\varepsilon (\alpha , \alpha) = -1 & \alpha \in \Phi\\
ii) &\varepsilon (\alpha , \beta) = - \varepsilon (\beta , \alpha) & \alpha,\beta,(\alpha+\beta)\in \Phi\qquad \text{antisymmetry}\\
iii) &\ep(\as ,\bb) = \ep(\bb, \as + \bb) & \text{if } \as, \as+\bb \in \Phi\, , \ \bb\in \Ll \\
iv) &\ep(\as ,\bb) = \ep(\bb, \as - \bb) & \text{if } \as, \as-\bb \in \Phi\, , \ \bb\in \Ll \\
\eeas
\label{remas}
\end{proposition}


\section{3-5 gradings and 3-linear products}\label{s:35}

\begin{figure}[h!]\centering
\includegraphics[width=60mm]{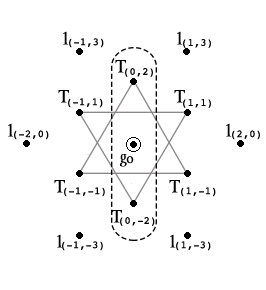}
	\caption{The Magic Star and a 3-graded subalgebra}\label{fig:mstg}
\end{figure}

\begin{figure}[h!]\centering
\includegraphics[width=50mm]{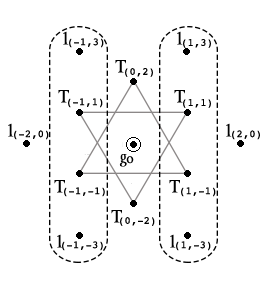}
	\caption{The Magic Star and a 5-grading}\label{fig:msfg}
\end{figure}

The Magic Star (from now on abbreviated by MS) reveals the presence of both a 3-graded and a 5-graded structure, see figures \ref{fig:mstg}, \ref{fig:msfg}. We state this in the following two Propositions, whose proof easily follows from the figures and the tables \ref{t:magices}, \ref{t:magiceo}, for the cases $\esn$ and $\eon$ (similarly for $\estn$).

\begin{proposition}
Let us denote by $\ggo$ the subalgebra of $\lms$ lying in the center of MS and by $T_{(r,s)}$ the subalgebra of $\lms$ spanned by the generators associated to the $(r,s)$-set of roots in the tables \ref{t:magices} and \ref{t:magiceo}, for $(r,s)=(\pm 1,\pm1),(0,\pm2)$. Notice that $T_{(r,s)}$ is also a  subalgebra with trivial product: $[T_{(r,s)},T_{(r,s)}]=0$. Then the subalgebra $\ggt :=\ggo\oplus \cc \oplus T_{(r,s)} \oplus T_{(-r,-s)}$ is three graded, $\cc$ representing any complex multiple of the Cartan associated to the root along the axis of the grading (the $(r,s)$ direction in the tables \ref{t:magices} $\div$ \ref{t:magiceo}). The subalgebra $\ggo\oplus \cc$ is in grade 0 and $T^+:=T_{(r,s)}$ ($T^-:=T_{(-r,-s)}$) is in grade 1 (-1). If $\lms=\eon$, in \cite{EP1} we noticed that $\ggt=\estn$.
\end{proposition}

\begin{table}[h!]
\begin{equation*}
\begin{array}{| c l c | c | c |}
\hline
\hspace{2.9cm} generalized\ roots & & (r,s) & \#\ of\ roots\\
\hline
\pm (k_1 -  k_2) &  & \pm (2,0) & 2\\
\pm (k_2 -  k_3) &  & \pm (-1,3) & 2\\
\pm (k_3 -  k_1) & & \pm (-1,-3) & 2\\
\hline
&&&\\
 \pm k_i \pm  k_j & 4 \le i < j \le N-3 &  & 2(N-6)(N-7)\\
&&(0,0)&\\
\frac{1}{2} (\pm (k_1 + k_2 + k_3) \pm k_4 \pm ... \pm {\bf u}) & \text{even \# of +} & & 2^{N-5} \\
&&&\\
\hline
&&&\\
k_1 +  k_2 \ , \  - k_3 \pm k_i & i=4,... , N-3 &  & 2N-11\\
&&(0,2)&\\
\frac{1}{2} (k_1 + k_2 - k_3 \pm k_4 \pm ... \pm {\bf u}) & \text{even \# of +} &  & 2^{N-6} \\
&&&\\
\hline
&&&\\
- k_1 -  k_2 \ , \  k_3 \pm k_i & i=4,... , N-3 &  & 2N-11\\
&&(0,-2)&\\
\frac{1}{2} (- k_1 - k_2 + k_3 \pm k_4 \pm ... \pm {\bf u}) & \text{even \# of +} & & 2^{N-6} \\
&&&\\
\hline
&&&\\
- k_2 -  k_3 \ , \  k_1 \pm k_i & i=4,... , N-3 &  & 2N-11\\
&&(1,1)&\\
\frac{1}{2} (k_1 - k_2 - k_3 \pm k_4 \pm ... \pm {\bf u}) & \text{even \# of +} & & 2^{N-6} \\
&&&\\
\hline
&&&\\
k_2 +  k_3 \ , \ - k_1 \pm k_i & i=4,... , N-3 &  & 2N-11\\
&&(-1,-1)&\\
\frac{1}{2} (- k_1 + k_2 + k_3 \pm k_4 \pm ... \pm {\bf u}) & \text{even \# of +} & & 2^{N-6} \\
&&&\\
\hline
&&&\\
- k_1 -  k_3 \ , \  k_2 \pm k_i & i=4,... , N-3 &  & 2N-11\\
&&(-1,1)&\\
\frac{1}{2} (- k_1 + k_2 - k_3 \pm k_4 \pm ... \pm {\bf u}) & \text{even \# of +} & & 2^{N-6} \\
&&&\\
\hline
&&&\\
k_1 +  k_3 \ , \  - k_2 \pm k_i & i=4,... , N-3 &  & 2N-11\\
&&(1,-1)&\\
\frac{1}{2} (k_1 - k_2 + k_3 \pm k_4 \pm ... \pm {\bf u}) & \text{even \# of +} & & 2^{N-6} \\
&&&\\
\hline
\end{array}
\end{equation*}
\vspace{-16pt}
\caption{The Magic Star for $\esn$; ${\bf u}:=k_{N-2}+k_{N-1}+k_N$}\label{t:magices}
\end{table}

\begin{proposition}
$\lms$ has a five-graded structure: take for instance $\ggt=\ggo\oplus \cc \oplus T_{(0,2)} \oplus T_{(0,-2)}$ then $\lms= (\ggt\oplus \cc)\oplus (T_{(1,1)} + 1_{(1,3)} + T_{(1,-1)} +1_{(1,-3)}) + (T_{(-1,1)} + 1_{(-1,3)} + T_{(-1,-1)} + 1_{(-1,-3)})\oplus 1_{(2,0)}\oplus 1_{(-2,0)}$,  where $\cc$ represents any complex multiple of the Cartan associated to the root $k_1-k_2$ associated to the grading and $1_{(r,s)}$ denotes the one-dimensional vector space spanned by $x_\as$, with $\as=\pm (k_2 -  k_3)$ for $(r,s)= \pm (-1,3)$; $\as=\pm (k_3 -  k_1)$ for $(r,s)= \pm (-1,-3)$ and $\as=\pm (k_1 -  k_2)$ for $(r,s)= \pm (2,0)$, see tables \ref{t:magices} $\div$ \ref{t:magiceo}.\\
The subalgebra $\ggt\oplus \cc$ is in grade 0; $F^+:=(T_{(-1,1)} + 1_{(-1,3)} + T_{(-1,-1)} + 1_{(-1,-3)})$ and $F^-:=(T_{(-1,1)} + 1_{(-1,3)} + T_{(-1,-1)} + 1_{(-1,-3)})$ are in grade 1 and -1; $1_{(\pm 2,0)}$ is in grade $\pm 2$.
\end{proposition}

\begin{table}[h!]
\begin{equation*}
\begin{array}{| c l c | c | c |}
\hline
\hspace{2.9cm} generalized\ roots & & (r,s) & \#\ of\ roots\\
\hline
\pm (k_1 -  k_2) &  & \pm (2,0) & 2\\
\pm (k_2 -  k_3) &  & \pm (-1,3) & 2\\
\pm (k_3 -  k_1) & & \pm (-1,-3) & 2\\
\hline
&&&\\
 \pm k_i \pm  k_j & 4 \le i < j \le N &  & 2(N-3)(N-4)\\
&&(0,0)&\\
\frac{1}{2} (\pm (k_1 + k_2 + k_3) \pm k_4 \pm ... \pm k_N) & \text{even \# of +} & & 2^{N-3} \\
&&&\\
\hline
&&&\\
k_1 +  k_2 \ , \  - k_3 \pm k_i & i=4,... , N &  & 2N-5\\
&&(0,2)&\\
\frac{1}{2} (k_1 + k_2 - k_3 \pm k_4 \pm ... \pm k_N) & \text{even \# of +} &  & 2^{N-4} \\
&&&\\
\hline
&&&\\
- k_1 -  k_2 \ , \  k_3 \pm k_i & i=4,... , N &  & 2N-5\\
&&(0,-2)&\\
\frac{1}{2} (- k_1 - k_2 + k_3 \pm k_4 \pm ... \pm k_N) & \text{even \# of +} & & 2^{N-4} \\
&&&\\
\hline
&&&\\
- k_2 -  k_3 \ , \  k_1 \pm k_i & i=4,... , N &  & 2N-5\\
&&(1,1)&\\
\frac{1}{2} (k_1 - k_2 - k_3 \pm k_4 \pm ... \pm k_N) & \text{even \# of +} &  & 2^{N-4} \\
&&&\\
\hline
&&&\\
k_2 +  k_3 \ , \ - k_1 \pm k_i & i=4,... , N &  & 2N-5\\
&&(-1,-1)&\\
\frac{1}{2} (- k_1 + k_2 + k_3 \pm k_4 \pm ... \pm k_N) & \text{even \# of +} & & 2^{N-4} \\
&&&\\
\hline
&&&\\
- k_1 -  k_3 \ , \  k_2 \pm k_i & i=4,... , N &  & 2N-5\\
&&(-1,1)&\\
\frac{1}{2} (- k_1 + k_2 - k_3 \pm k_4 \pm ... \pm k_N) & \text{even \# of +} &  & 2^{N-4} \\
&&&\\
\hline
&&&\\
k_1 +  k_3 \ , \  - k_2 \pm k_i & i=4,... , N &  & 2N-5\\
&&(1,-1)&\\
\frac{1}{2} (k_1 - k_2 + k_3 \pm k_4 \pm ... \pm k_N) & \text{even \# of +} & & 2^{N-4} \\
&&&\\
\hline
\end{array}
\end{equation*}
\vspace{-16pt}
\caption{The Magic Star for $\eon$}\label{t:magiceo}
\end{table}

\begin{definition}
The 3-grading $\ggt :=\ggo\oplus \cc \oplus T^+ \oplus T^-$allows us to introduce a 3-linear product $(\ ,\ ,\ ):T^+\times T^+\times T^+\to T^+$, once we define a conjugation $\zeta: \tp \to \tm$:
\be
(x,y,z) := [[x,\zeta(y)],z]
\ee
alternatively we can define two 3-linear products on the pair $(T^+,T^-)$:
\be
(x^\sigma,y^{-\sigma},z^\sigma) := [[x^\sigma,y^{-\sigma}],z^\sigma]\ ,  \ x^\sigma, z^\sigma \in T^\sigma \, , \, y^{-\sigma}\in T^{-\sigma} \ , \ \ \sigma =\pm
\ee
Analogously the 5-grading and a conjugation $\kappa:F^+\to F^-$, allow us to introduce a 3-linear product  $(\ ,\ ,\ ):F^+\times F^+\times F^+\to F^+$:
\be
(x,y,z):=[[x,\kappa (y)],z]
\ee
or two 3-linear products on the pair $(F^+,F^-)$:
\be
(x^\sigma,y^{-\sigma},z^\sigma) := [[x^\sigma,y^{-\sigma}],z^\sigma]\ ,  \ x^\sigma, z^\sigma \in F^+ \, , \, y^{-\sigma}\in F^- \ , \ \ \sigma =\pm
\ee
\end{definition}

\begin{remark}
We remark that the spaces in grade $\pm 2$ are one-dimensional, namely the 5-grading is of contact type \cite{Santi}.
\end{remark}
\begin{definition}
Given a 5-graded Lie algebra $\gmf = \gmf_{-2} + \gmf_{-1} + \gmf_0 + \gmf_1 + \gmf_2$ over a field $\fff$ with one dimensional $\gmf_{\pm2} := \fff\, x_{\pm \rho}$ there is a natural way, \cite{hel} to define a non-degenerate skew-symmetric bilinear $<\ , \ >$ form and a fully symmetric quartic form $q$ on $\gmf_1$ by:
\bea{ll}
[ x, y]= <x,y> x_\rho & x,y\in \gmf_1\\
q(x,x,x,x)= [x , [x , [x , [x , x_{-\rho}]]]] &  x\in \gmf_1\\
q(x_{\beta_1}, x_{\beta_2} , x_{\beta_3} , x_{\beta_4})= \dfrac{1}{4!} \sum_{\pi \in S_4}{
[x_{\beta_{\pi (1)}} , [x_{\beta_{\pi (2)}} , [x_{\beta_{\pi (3)}} , [x_{\beta_{\pi (4)}} , [ x_{-\rho}]]]]} &x_{\beta_i}\in \gmf_1
\eea
where $S_4$ is the symmetric group in $\{ 1,2,3,4 \}$.
\end{definition}

\begin{remark}
For $n=1$ the subalgebras $T^+,T^-$ are the Jordan algebras of $3\times 3$ Hermitian matrices over the complex numbers in the case of $\es$, the quaternions in the case of $\est$, and the octonions in the case of $\eo$. With the 3-linear product $T^+$ becomes a Jordan triple system,  \cite{loos1}.The pair $(T^+,T^-)$ is a Jordan Pair,  \cite{loos1}. Still for $n=1$, $F^+$ is a Freudenthal triple system, \cite{hel} and $(F^+,F^-)$ is a Kantor Pair, \cite{af1}.
\end{remark}

\section{{\em HT}-algebras}\label{s:dta}

We now concentrate on the set of roots $T_{(r,s)}$ on a tip of the Magic Star. In order to fix one, let us consider $T_{(1,1)}$ and denote it simply by {\em T}. The reader will forgive us if we shall be concise whenever there is no risk of misunderstanding and use $T$ to denote both the set of roots and the set of elements in $\lms$ associated to those roots. An element of $T$ is an $\fff$-linear combination of $x_\as, x_\beta, ...$ for $\as,\bb,...$ in $T_{(1,1)}$. We give $T$ an algebraic structure with a commutative product, thus mimicking the case $n=1$ when $T$ is a Jordan algebra.

\begin{figure}[h]
\centering
\includegraphics[width=0.60\textwidth]{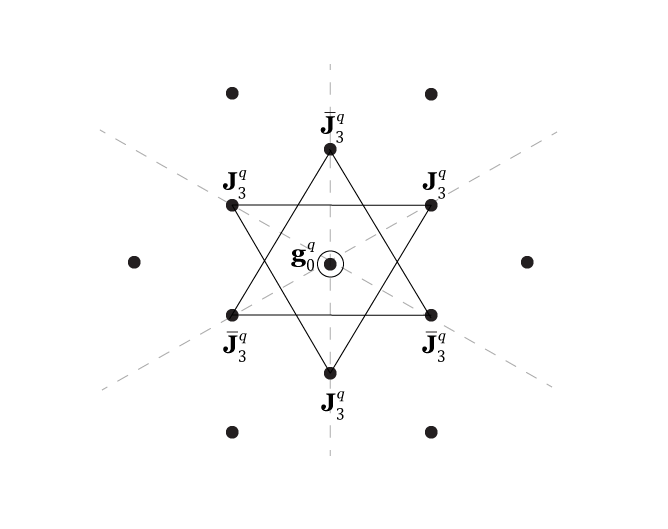}
\caption{The Magic Star of finite dimensional
exceptional Lie algebras \protect\cite{Mukai, Truini}. $\mathbf{J}_{3}^{q}$
denotes a simple Jordan algebra of rank-3, parametrized by $q=\dim _{\mathbb{%
R}}\mathbb{A}=1,2,4,8$ for $\mathbb{A}=\mathbb{R},\mathbb{C},\mathbb{H},%
\mathbb{O}$, corresponding to $\mathbf{f}_{4}$, $\mathbf{e}_{6}$, $\mathbf{e}%
_{7}$, $\mathbf{e}_{8}$, respectively. In the case of $\mathbf{g}_{2}$
(corresponding to $q=-2/3$), the Jordan algebra is trivially the identity
element : $\mathbf{J}_{3}^{-2/3}\equiv \mathbb{I}:=diag(1,1,1)$.}
\label{fig:MagicStar}
\end{figure}

\begin{remark}
Let us first show what happens in the case $n=1$, in particular for $\eo$. It is proven in \cite{Marrani-Truini-1} that $T=J_{(1,1)}$ is the Jordan Algebra $\joto$ of $3\times 3$ Hermitian matrices over the octonions. Let us denote by
$P_1 := E_{11}$, $P_2 := E_{22}$, $P_3 := E_{33}$ the three trace-one idempotents, whose sum is the identity in $\joto$, $E_{ii}$ representing the matrix with a 1 in the $(ii)$ position and zero elsewhere.\\
Let us identify the algebra $\dq \subset \fq\subset\es\subset\es\oplus\adc\subset\eo$ with the roots $\pm k_i \pm k_j$, $4\le i<j \le 7$ and  $P_1,P_2,P_3$ with the elements of $J_{(1,1)}$ which are left invariant by $\mathbf{d_4}$, \cite{jacob2}. This uniquely identifies $P_1 , P_2 , P_3$ with the roots $k_1+ k_8$, $k_1- k_8$ and $-k_2-k_3$. Since $\est\subset\eo$ is 3-graded, one can define a 3-linear product on $J_{(1,1)}$, \cite{mac1}, and a Jordan product through the correspondence with the quadratic formulation of Jordan algebras.

Still in the case of $\eo$ ($\eon$ for $n=1$), if we identify the $\eo$ generators corresponding to the roots in $\Phi_O$ with the bosons and those corresponding to the roots in $\Phi_S$ with the fermions, then the 11 bosons in $T$ are $k_1 \pm k_j$, $j=4,...,8$ and $-k_2-k_3$. The three bosons $k_1\pm k_8$ and $-k_2-k_3$, correspond to the three diagonal idempotents $\xu,\xd,\xt$, which are left invariant by $\dq$, whose roots are $\pm k_i \pm k_j \ , \ 4\le i<j \le 7$.\\
The remaining 8 vectors (bosons) and the 2x8 spinors (fermions) are linked by triality. Notice that one spinor has $\frac12(k_1-k_2-k_3 + k_8)$ fixed and even number of + signs among $\pm k_4 \pm k_5\pm k_6\pm k_7$, while the other one has $\frac12(k_1-k_2-k_3 - k_8)$ fixed and odd number of + signs among $\pm k_4\pm k_5\pm k_6\pm k_7$. So they are both 8 dimensional representations of $\dq$.\\
Therefore a generic element of $\joto$ is a linear superposition of 3 diagonal elements plus one vector and two spinors (or a bispinor). The vector can be viewed as the $(12)$ octonionic entry of the matrix (plus its octonionic conjugate in the $(21)$ position) and the two 8-dimensional spinors as the $(31)$ and $(23)$ entry (plus their respective octonionic conjugate in the $(31)$ and $(23)$ position).

\textit{Mutatis mutandis}, analogous results hold for all finite dimensional exceptional Lie algebra; the Magic Star projection (depicted in Fig. \ref{fig:MagicStar}) of the corresponding
root lattice onto a plane determined by an $\mathbf{a}_{\mathbf{2}}$ root
sub-lattice has been introduced by Mukai \cite{Mukai}, and later investigated in depth in \cite{Truini}
(see also \cite{Marrani-Truini-1}), with a different approach exploiting Jordan Pairs
\cite{loos1}; see also the discussion in \cite{Group32-1} and \cite{Group32-2}.
\end{remark}

One can then do the same for $n\ge 1$ and consider the most general case $\eon$ (the other two cases $\esn$ and $\estn$ being a restriction of it; for a detailed treatment of $\fqn$, see \cite{EP3}).

Let us denote $\xu$, $\xd$ and $\xt$ the elements of $\lms$ in $T$ associated to the roots $\rho_1:=k_1+ k_N$, $\rho_2:=k_1- k_N$ and $\rho_3:=-k_2-k_3$:
\be\label{notp}
\xu \lra \rho_1:=k_1+ k_N \ ; \ \xd \lra \rho_2:=k_1- k_N \ ; \ \xt \lra \rho_3:=-k_2 - k_3
\ee
They are left invariant by the Lie subalgebra $\mathbf{d_{N-4}}=\dqn$, whose roots are $\pm k_i \pm k_j \ , \ 4\le i<j \le N-1$.\\
We denote by $\tfo$ the set of roots in $T\cap \Phi_O$,  by $\tfop$ the set of roots in $\tfo$ that are not $\rho_1,\rho_2,\rho_3$ and by $\tfs$ the set of roots in $T\cap \Phi_S$. In the case we are considering, where $T=T_{(1,1)}$ we have $\tfop = \{ k_1\pm k_j\, ,\ j=4,...,N-1\}$ and $\tfs= \{ \frac{1}{2} (k_1 - k_2 - k_3 \pm k_4 \pm ... \pm k_N)\}$, even $\#$ of $+$. We further split $\tfs$ into $\tfsp = \{ \frac{1}{2} (k_1 - k_2 - k_3 \pm k_4 \pm ... + k_N)\}$ and $\tfsm= \{ \frac{1}{2} (k_1 - k_2 - k_3 \pm k_4 \pm ... - k_N)\}$. Then $v\in \tfop$ is an $8n$-dimensional vector and $s^\pm$ are $2^{4n-1}$-dimensional chiral spinors of $\dqn$.

We write a generic element $x$ of $T$ as  $x= \sum \lambda_i \xpi + \, x_v\, +\, x_{s^+}+x_{s^-}$ where
\be
x_v = \sum_{\as \in \tfop} \lambda_\as^v x_\as
\ee
\be
x_{s^\pm} = \sum_{\as \in \tfspm} \lambda_\as^{s^\pm} x_\as
\ee

We view $\lambda_\as^v$ as a coordinate of the vector $\lambda^v$ and $\lambda_\as^{s^\pm}$ as a coordinate of the spinor $\lambda^{s^\pm}$. The conjugation properties of the vector and chiral spinor representations of $\dqn$ imply the existence of the symmetric bilinear forms $\omv$ and $\oms$ (whose explicit expressions are computed in \cite{EP4}). We will denote by  $\bar\lambda^v$ ($\bar\lambda^{s^\pm}$) the vector (spinor) in the dual space with respect to such bilinear forms. We can therefore view $x$ as a $3\times 3$ Hermitian matrix:

\be\label{matrix}
\left(
\begin{array}{lll}
\lambda_1 &\lambda_v &\bar  \lambda_{s^+}\\
\bar \lambda_v &\lambda_2 &\lambda_{s^-}\\
\lambda_{s^+} &\bar  \lambda_{s^-} &\lambda_3
\end{array}
\right)
\ee
 whose entries have the following $\fff$-dimensions:
\bit
\item 1 for the {\it scalar} diagonal elements $\lambda_1,\lambda_2,\lambda_3$;
\item $8n$ for the {\it vector} $\lambda_v$;
\item $2^{4n-1}$ for the {\it spinors} $\lambda_{s^\pm}$.
\eit

\begin{remark}
We see that only for $n=1$ the dimension of the vector and the spinors is the same, whereas for $n>1$ the entries in the $(12), (21)$ position have different dimension than those in the $(31),(13),(23),(32)$ position. Nevertheless, we define in equation \eqref{jprod} a commutative product of the elements in $T$, which then becomes a generalization of the Jordan algebra $\joto$ in a very precise sense, which we will henceforth name {\em HT}-algebra. The set of matrices \eqref{matrix} is contained in a cubic matrix algebra \cite{Vinberg}, whose product is the standard block matrix product; blockwise, by the representation theory of $\dqn$, such a product suitably exploits the gamma matrices (also named Dirac matrices) (see e.g. \cite{gamma}) of $\dqn$ itself. In a forthcoming paper \cite{EP4}, we will investigate this matrix algebra further; we think that this is indeed a rank-3 {\em T}-algebra \cite{Vinberg}.
\end{remark}

\begin{figure}[h]
\centering
\includegraphics[width=0.60\textwidth]{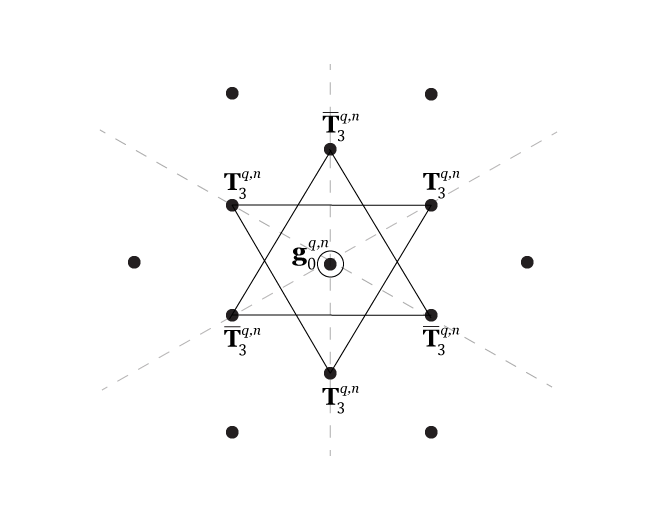}
\caption{The Magic Star projection Magic Star algebras. $\mathbf{T}%
_{3}^{q,n}$ denotes an {\em HT}-algebra \cite{Vinberg, EP4}, parametrized by $q=\dim _{\mathbb{R}}\mathbb{A}%
=1,2,4,8$ for $\mathbb{A}=\mathbb{R},\mathbb{C},\mathbb{H},\mathbb{O}$, and $%
n\in \mathbb{N}$, corresponding to $\mathbf{f}_{4}^{(n)}$ \cite{EP3}, $%
\mathbf{e}_{6}^{(n)}$, $\mathbf{e}_{7}^{(n)}$, $\mathbf{e}_{8}^{(n)}$,
respectively. We remind that the smallest exceptional Lie algebra $\mathbf{g}_{2}$
(corresponding to $q=-2/3$) cannot be EP-generalized, because it does not
enjoy a spin factor embedding \cite{Group32-1}, and because $\mathbf{J}_{3}^{-2/3}\equiv
\mathbb{I}$ \cite{Mile-High-EP, EP1}.}
\label{fig:MagicStarT}
\end{figure}

\begin{remark}
Considering all non-trivial Magic Star algebras \cite{EP1} (i.e. $\fqn$, $\esn$, $\estn$ and $\eon$),  we can state the following result (anticipated in \cite{Mile-High-EP}):
there exists an\footnote{%
Such a projection is not unique; as for the Magic Star projection of exceptional Lie algebras,
depicted in Figure \ref{fig:MagicStar}, also for the Magic Star projection of Magic Star
algebras, depicted in Figure \ref{fig:MagicStarT}, there are four possible
equivalent projections.} $\mathbf{a}_{2}$ projection of the Magic Star algebras
(namely, of the corresponding generalized root lattice), such that \textit{a Magic Star
structure persists}, with suitable generalizations of rank-$3$ Jordan
algebras $\jotqq$ (given by the aforementioned {\em HT}-algebras) occurring on the tips of the persisting Magic Star! The
 Magic Star projection of Magic Star algebras is depicted in Fig. \ref{fig:MagicStarT}.
\end{remark}

\begin{definition}
We denote by $I$ the element $I:=\xu+\xd+\xt$ and by $I^-$ the {\it Cartan involution} of $I$, namely the element $I^-:=-\xub-\xdb-\xtb$ of $\bar T:= T_{(-1,-1)}$, where $\xub$, $\xdb$ and $\xtb$ are associated to the roots $-k_1- k_N$, $-k_1+ k_N$ and $k_2+k_3$.\\

We give $T$ an algebraic structure by introducing the commutative product, see Proposition \sref{t-al}:
\be\label{jprod}
x\jp y := \frac12 [[x,I^-],y] \quad ,\quad x,y\in T
\ee

We introduce the trace (see Proposition \sref{t-al}) $tr(x)\in \fff$ for $x\in T$ in the following way:
\be\label{trace}
\text{for } x=\ell_1 \xu+\ell_2 \xd+\ell_3 \xt +\sum_{\substack{\as\ne\rho_1,\rho_2,\rho_3}}\ell_\as x_\as\quad tr(x)=\ell_1+\ell_2+\ell_3
\ee

We denote by $tr(x,y):=tr(x\jp y)$ and by $x^2:=x\jp x$. For each $x\in T$ we define
\be \label{sharp}
x^\# = x^2 -tr(x) x - \frac12 (tr(x^2) - tr(x)^2) I
\ee
and we say that $x$ is rank-1 if $x^\#=0$.

Let us also introduce $N(x)\in \fff$ for $x\in T$ in the following way:
\be\label{norm}
N(x) = \frac16 \left\{ tr(x)^3 - 3 tr(x)tr(x^2) +2 tr(x^3) \right\} = \frac13 tr(x^\#,x)
\ee
where $x^2= x\jp x$ and $x^3= x\jp x^2$.
\end{definition}

\begin{remark}
Notice that $tr(x^\#)=- \frac12 (tr(x^2) - tr(x)^2)$, therefore $x^\#=0$ implies $x^2 = tr(x) x$: a rank-1 element of $T$ is either a nilpotent or a multiple of a {\it primitive} idempotent $u\in T$: $u^2 = u$ and $tr(u)=1$.
\end{remark}

\begin{remark}
Since by \eqref{norm} a rank-1 element $x$ has $N(x)=0$ any element of $T$ falls into the following classification:
\bea{ll}\label{class}
\text{rank-1:} & x^\# = 0\\
\text{rank-2:} & x^\# \ne 0 \ , \ N(x)=0\\
\text{rank-3:} & N(x)\ne 0\\
\eea

Examples of rank-2 and rank-3 elements are $\xu + \xd$ and $I$ respectively.
\end{remark}

\begin{definition}
Let us define $\dd$ the $\lms$ Lie subalgebra $\ddnt$ if $\lms=\esn$, $\ddnd\oplus\au$ if $\lms=\estn$, $\ddn$ if $\lms=\eon$, spanned by $\{x_\alpha \, , \ \alpha \in \Phi_O\}$; see remark $\mathbf{3.1}$ of \cite{EP1}.
\end{definition}

We can prove the following

\begin{proposition}\label{t-al}
The product $x,y\to x\jp y$ is commutative. With respect to this product the element $I\in T$ is the identity  and the elements $\xu , \xd , \xt\in T$ are trace-one idempotents, hence are rank-1 elements of $T$. All $x_\as \in T$ but  $\xu , \xd , \xt$ are nilpotent and traceless, hence they are also rank-1 elements of $T$. The form $x \to tr(x)$ is a trace form, namely $tr(x,y)$ is bilinear and symmetric in $x,y$ and the associative property $tr(x\jp y,z)=tr(x,y\jp z)$ holds for every $x,y,z \in T$.
\end{proposition}
\proof If $x,y\in T$ they are linear combinations of elements $x_\alpha, x_\beta,...$ associated to roots $\alpha,\beta,...$ whose sum is not a root. Therefore $[x,y]=0$. Moreover, since $I^-\in \dd$, Proposition 6.1 of \cite{EP1} yields that $\adj_{I^-}$ is a derivation, hence:
\be [[x,I^-],y] + [[I^-,y],x] = 0\ \Rightarrow\ [[x,I^-],y] = [[y,I^-],x]
\ \Rightarrow\ x\jp y = y\jp x
\ee
We also have
\be
[\xpi,-\xpjb]=\delta_{ij} \hri \ , \ i,j=1,2,3
\ee
and for any $\as \in T_{(1,1)}$, being $\rho_1+\rho_2+\rho_3 = 2k_1-k_2-k_3$, see tables \ref{t:magices} and \ref{t:magiceo}:
\be
I\jp x_\as = \frac12 [[I,I^-],x_\as] = \frac12 [\hru+\hrd+\hrt,x_\as] = \frac12 (\rho_1+\rho_2+\rho_3,\as) x_\as = x_\as
\ee
By linearity this extends to any $x\in T$.\\
Moreover:
\be \xpi\jp \xpi=\frac12 [[\xpi,I^-],\xpi] = \frac12 [h_{\rho_i},\xpi] = \frac12 (\rho_i,\rho_i) \xpi = \xpi
\ee
Obviously $tr(\xu)=tr(\xd)=tr(\xt)=1$.\\

We now show that $x_\as \in T$ is nilpotent if $\as\ne k_1\pm k_N, -k_2-k_3$.\\
Suppose $x_\as\jp x_\as\ne 0$ then either $2 \as - k_1 - k_N$ or $2 \as - k_1 + k_N$ or $2 \as + k_2 + k_3$ is a root. If $\as\in \Phi_O$ then $2 \as - k_1 \mp k_N = k_1 \pm 2 k_i \mp k_N$ is a root if and only if $k_i=k_N$, whereas $2\as + k_2 +k_3$ is not a root. If $\as\in \Phi_S$ then nor $2 \as - k_1 \mp k_N$ nor $2\as + k_2 +k_3$ are roots. Therefore $\as\ne k_1\pm k_N, -k_2-k_3$ implies $x_\as^2=0$ and, obviously, $tr(x_\as)=0$ and rank($x_\as$)$=1$.\\

Finally, from the definition of the product $x,y\to x\jp y$ and the definition of the trace it easily follows that $tr(x,y)$ is bilinear. The associativity of the trace is proven in Proposition \sref{assoc-tr} in App. \sref{s:astr}. \esquare \\

\section{{\em HT}-pairs}\label{s:pairs}
\begin{definition}
We define {\em HT}-pair the pair $\pt =(T^+,T^-)$ where $T^+$ ($T^-$) is the $\fff$-vector space generated by $\{x^+_\as \ | \as\in T_{(1,1)}$ a(resp. $\{x^-_\as \ | \as\in T_{(-1,-1)}$) so that the roots associated to $T^+$ and $T^-$ are on opposite tips of the Magic Star.
Define the maps $U_{x^\pm}y^\mp = \frac12[[x^+,y^-],x^+]$ quadratic in $x^\pm$ and linear in $y^\mp$, such that $U_{x^\pm}:T^\mp \to T^\pm$. Define the linearization $V$ of $U$: $V_{x,y}z = (U_{x+z}-U_x -U_z)y = \frac12([[x,y],z]+[[z,y],x])$ for $x,z\in T^\pm$ and $y\in T^\mp$. We have $V_{x,y}x= 2 U_x y$.
\end{definition}
\begin{remark}
We remark that only in the case $n=1$, due to the Jacobi identity, $V_{x,y}z = [[x,y],z]$ for all $x,z\in T^\pm$, $y\in T^\mp$. In this case, the {\em HT}-pairs become Jordan pairs \cite{loos1}.
\end{remark}
\begin{definition}
We define an {\em idempotent} of the {\em HT}-pair $\pt$ as the pair $(x,y)$ such that $U_x y = x$ and $U_y x=y$.
\end{definition}

Then, one can prove the following
\begin{proposition}
The rank-1 elements $x_\as$ for $\as \in T$, be they nilpotent or primitive idempotent, can be completed to become idempotents of $\pt$.
\end{proposition}
\proof For $x=x_\as$ take $y=-\dfrac2{(\as,\as)}x_{-\as}$ then $U_xy=\dfrac1{(\as,\as)}[h_\as,x_\as]=x_\as$ and $U_yx=\dfrac1{(\as,\as)}\dfrac2{(\as,\as)}[h_{\as},x_{-\as}]=-\dfrac2{(\as,\as)}x_{-\as}=y$. \esquare \\

For many aspects the language of {\em HT}-pairs inside $\eon$ is more natural than that of {\em HT}-algebras. A detailed investigation of the properties of such structures is beyond the scope of the present paper, and it is left for further future work.


\section{Further Developments}\label{s:outlook}

A number of research venues stems from the novel non-Lie, infinitely
countable chains of finite dimensional generalizations of exceptional Lie
algebras provided by Magic Star algebras, introduced in \cite{EP1} and in
the present paper, after the anticipation given in \cite{Mile-High-EP}.
We here list some of the developments which we plan to report on in the
near future.

Firstly, we will consider in detail the violation of the Jacobi identity
determining the non-Lie nature of Magic Star algebras. As anticipated in
\cite{EP1}, such a violation occurs due to the non-trivial (i.e.
non-Abelian) nature of the spinorial subsector of the Magic Star algebras.
In fact, Magic Star algebras are not simply non-reductive extensions of
(pseudo-)orthogonal Lie algebras by some translational (i.e. Abelian) spinor
generators (as it holds for the electric-magnetic duality Lie algebra of $%
\mathcal{N}=2$ supergravity coupled to homogeneous non-symmetric scalar
manifolds in $D=3,4,5$ Lorentzian space-time dimensions \cite{dWVP,dWVVP}),
but rather they are characterized by non-translational (i.e., non-Abelian)
spinorial extensions, whose non-trivial commutation relations with all
generators generalize, in the very spirit of Exceptional Periodicity, the
commutation relations characterizing exceptional Lie algebras in their
\textquotedblleft spin factor embedding decompositions" (cfr. \cite%
{Group32-1}, as well as Remark 3.1 of \cite{EP1}).

As pointed out also in \cite{EP1}, we remind that Exceptional Periodicity
provides a way to go beyond $\mathbf{e}_{8}$ (and the whole set of
exceptional Lie algebras but $\mathbf{g}_{2}$) which is very different from
the way provided by affine and (extended) Kac-Moody algebras \cite{kac},
which also appeared as symmetries for (super) gravity models reduced to $%
D=2,1,0$ dimensions (see e.g. \cite{extended-Refs, West}), as well as near
spacelike singularities in supergravity \cite{spacelike-Refs}. In the near
future, we plan to analyze in detail the relation between Magic Star
algebras and (generalized) Kac-Moody (and possibly Borcherds) algebras.

Also, in the forthcoming paper \cite{EP3} we will deal with the
determination of the inner derivations of {\em HT}-algebras, highlighting the relation with the Magic Star algebra $\mathbf{f}%
_{4}^{(n)}$; subsequently, we will study further larger symmetry algebras
(such reduced structure, conformal and possibly quasi-conformal symmetries),
analyzing similarities and differences with the electric-magnetic duality
Lie algebras of $\mathcal{N}=2$ supergravity coupled to homogeneous
non-symmetric scalar manifolds in $D=3,4,5$ Lorentzian space-time dimensions
\cite{dWVP,dWVVP}.

Concerning further physical applications, we envisage potential applications
of Exceptional Periodicity and Magic Star algebras in Quantum Gravity,
stemming from the considerations in \cite{Truini}, \cite%
{Marrani-Truini-Interactions, Group32-2, Truini-EUG}. We aim at a consistent model of
elementary particle physics and space-time expansion at the early stages of
the Universe based on the idea that interactions, defined in a purely
algebraic way, are the fundamental objects of the theory, whereas
space-time, hence gravity, is an emergent entity. Within this perspective,
the failure of the Jacobi identity in the spinor sector of Magic Star
algebras might be related to dark matter /dark energy degrees of freedom,
and the non-Abelian nature of the spinor component of Magic Star algebras
would play a key role in order to have non-trivial interactions among bosons
and fermions in such algebras, which are not superalgebras nor $\mathbb{Z}%
_{2}$-graded algebras.

In our view, a synthesis of spectral algebraic geometry (suitably
generalized in a non-associative sense within {\em HT}-algebras) and
Exceptional Periodicity can yield to the unification of various approaches
to Quantum Gravity, thus providing a new lens for searching a
non-perturbative theory of all matter, forces and space-time.

\section*{Acknowledgements}
We thank Willem de Graaf for useful discussions on {\em T}-algebras and enlightening feedback on an early version of the manuscript.

\appendix

\section{Associativity of the trace}\label{s:astr}

We recall that in section \sref{s:dta} we have fixed $T:=T_{(1,1)}$. Therefore we refer to the row $(1,1)$ of table \ref{t:magiceo} for the explicit expression of the roots involved in T and $-k_1\pm k_N$, $k_2+k_3$ for the roots in $I^-$ appearing in the commutative product \eqref{jprod}.\\
In the proofs that will follow, we make use of the following identities (see (3.10), (3.11) and (3.13) of \cite{EP1}, here with $N=R$):
\bea{ccl} \label{kal}
k_{N-1} &= &\frac12 (\alpha_{N-1} - \alpha_{N-2})\\ \\
k_i &= &\alpha_i + k_{i+1} = \sum_{\ell=i}^{N-2}\alpha_\ell +\frac12 (\alpha_{N-1}-\alpha_{N-2}) \ , \ 1\le i\le N-2\\ \\
{\bf k_N} &= &- 2 \alpha_N - \sum_{\ell=1}^{N-2} {\ell \alpha_\ell - \frac{N-1}{2}(\alpha_{N-1}}-\alpha_{N-2})
\eea
from which we obtain, for $1\le i<j\le N-1$ and forcing $\sum_{\ell=r}^s{\alpha_\ell}=0$ if $r>s$:

\be\label{posroots1}
\left.
\begin{array}{rcl}
k_i + k_j &= & \sum_{\ell=i}^{N-3}\alpha_\ell + \sum_{\ell=j}^{N-2}\alpha_\ell + \alpha_{N-1}\\ \\
k_i - k_j &= & \sum_{\ell=i}^{j-1}\alpha_\ell
\end{array}
\right\}
1\le i<j\le N-1
\ee
\bea{lcl}
\pm k_i - k_N &= & 2 \alpha_N +\sum_{\ell=1}^{i-1} {\ell \alpha_\ell}+ \sum_{\ell=i}^{N-3} {(\ell\pm1)\alpha_\ell} + (2n+ \frac{1 \pm 1}{2})\alpha_{N-2}\\ \\
&&+ (2n+1 + \frac{1 \pm 1}{2})\alpha_{N-1} \ , \quad i\le N-2\\ \\
\pm k_{N-1} - k_N &= & 2 \alpha_N + \sum_{\ell=1}^{N-3} {\ell \alpha_\ell} + (2n+\frac{1 \mp 1}{2})\alpha_{N-2}+ (2n+1 + \frac{1 \pm 1}{2})\alpha_{N-1}
\label{posroots3}
\eea

Let us introduce the notation $x^\pm_{v\mu}$, $\mu=0,...,N-5$, to denote the vectors:
\be\label{notv}
x^\pm_{v\mu} = \left\{
\ba{ll}
x^\pm_{v i}:=x_{k_1\pm k_{i+3}} & 1\le i\le N-5 \\
x^\pm_{v o}:= \pm\, x_{k_1\pm k_{N-1}}
\ea
\right.
\ee
A generic vector in $T$ is a linear combination of the vectors $x^\pm_{v\mu}$. We need the following

\begin{lemma}\label{xpv}
We have the following products among scalars and vectors in the notations \eqref{notp}, \eqref{notv}:
\bea{llcl}
i) & \xpi\jp \xpj &=& \delta_{ij} \xpi \\
ii) & \xpi\jp x_\bb &=& \frac12 (\rho_i,\bb) x_\bb \\
iii) & x^\pm_{v\mu}\jp x^\mp_{v\mu} &=&\frac12 (\xu+\xd) \\
iv) & x^\pm_{v\mu}\jp x^\mp_{v\nu} &=& 0 \quad \mu\ne \nu \\
v) & x^\pm_{v\mu}\jp  x^\pm_{v\nu} &=& 0 \\
\eea
\end{lemma}
\proof
The first two equations $i)$ and $ii)$ are easily proven:
\bea{lcl}
\xpi\jp \xpj &=& \frac12 [[\xpi, -\xub-\xdb-\xtb], \xpj] = \frac12 [h_{\rho_i}, \xpj]= \delta_{ij} \xpi \\
\xpi\jp x_\bb &=& \frac12 [h_{\rho_i}, x_\bb] = \frac12 (\rho_i,\bb) x_\bb
\eea
We now prove $iii)$. Since $[x_{vi},\xtb]=0$ we get:
\bea{lcl}
x^\pm_{v\mu}\jp x^\mp_{v\mu} &=& (-1)^{\delta_{\mu0}}\frac12\left[-\ep(k_1\pm k_i,-k_1-k_N)\ep(\pm k_i - k_N,k_1 \mp k_i) \xd\right.\\ &&\left. -\ep(k_1 \pm k_i,-k_1+k_N)\ep(\pm k_i+k_N,k_1 \mp k_i) \xu \right]
\eea
We now show that:
\bea{l}\label{vvb}
\ep(k_1+k_i,-k_1\mp k_N)\ep(k_i \mp k_N,k_1- k_i)=\\
\ep(k_1-k_i,-k_1\mp k_N)\ep(-k_i \mp k_N,k_1+ k_i)=
\left\{ \begin{array}{l} -1 \text{ for } i\le N-2\\ \phantom{-} 1\text{ for } i = N-1\end{array} \right.
\eea

Using $i),ii),vi)$ of Proposition \sref{epsprop}, as well as Proposition \sref{remas}, being the arguments of both terms roots whose sum is a root, we get
\bea{l}\label{vpm1}
\ep(k_1+k_i,-k_1\mp k_N)\ep(k_i \mp k_N,k_1- k_i)=\\
\ep(k_1-k_i, k_i\mp k_N)\ep(- k_1 \mp k_N,k_1+ k_i)=\\
\ep(k_1-k_i, k_i +k_1 -k_1\mp k_N)\ep(- k_1 +k_i -k_i \mp k_N,k_1+ k_i)=\\
\ep(k_1-k_i,-k_1\mp k_N)\ep^2(k_1-k_i,k_1+k_i)\ep(-k_i \mp k_N,k_1+ k_i)= \\
\ep(k_1-k_i,-k_1\mp k_N)\ep(-k_i \mp k_N,k_1+ k_i)\\
\eea

Since $2k_i$ is in the lattice $\Ll$ we can add and subtract $2k_i$ in any argument of $\ep$. By doing so we obtain
\bea{l}\label{vpm2}
\ep(k_1-k_i,-k_1\mp k_N)\ep(-k_i \mp k_N,k_1+ k_i)= \\
-\ep(k_1-k_i,-k_1\mp k_N)\ep(k_1- k_i, k_i \pm k_N)\ep(2k_i,k_i \pm k_N)= \\
-\ep(k_1-k_i,k_i-k_1)\ep(2k_i,k_i \pm k_N)=\ep(2k_i,k_i \pm k_N)
\eea

Because of \eqref{epsdef}, whenever one argument of $\ep$ is an even multiple of a simple root, that argument can be neglected, as it would contribute with a factor $+1$. Therefore, by \eqref{kal} and \eqref{posroots3}, one obtains
\be\label{vpm3}
\ep(2k_i,k_i \pm k_N) = \ep(\as_{N-1}-\as_{N-2},k_i \pm k_N)
\ee

Since the ordering of the simple roots is $\alpha_i<\alpha_j$ if $i<j$, and since $\ep(\as_j,\as_i)=1$ if $\as_j>\as_i$ (hence $j>i$), $\ep(\as_i,\as_i)=-1$ and $\ep(\as_{N-2},\as_{N-1})=1$, because $\as_{N-2}+\as_{N-1}\notin \Phi$, we get, for $i\le N-2$
\be\label{vpm4}
\ep(2k_i,k_i \pm k_N) = \ep(\as_{N-1}-\as_{N-2},\dfrac{1\mp 1}2\as_{N-2}+\dfrac{3\mp 1}2\as_{N-1}) =-1 \, , \, i\le N-2
\ee
and for $i = N-1$
\be\label{vpm5}
\ep(2k_{N-1},k_{N-1} \pm k_N) = \ep(\as_{N-1}-\as_{N-2},\dfrac{1\pm 1}2\as_{N-2}+\dfrac{3\mp 1}2\as_{N-1}) =1 \, ,\, i = N-1
\ee
We have thus proven \eqref{vvb}, from which $iii)$ follows.\\

The equations $iv)$ and $v)$ are trivial, since, for $4\le i,j \le N-1$, $k_1\pm k_i-k_1\pm k_N + k_1 \pm k_j = k_1\pm k_i\pm k_j \pm k_N$ is not a root unless $k_j=-k_i$.\esquare\\

We can proceed to prove the main result of this Appendix.

\begin{proposition}\label{assoc-tr}
The form $x \to tr(x)$ is associative, namely $tr(x\jp y,z)=tr(x,y\jp z)$ holds for every $x,y,z \in T$.
\end{proposition}
\proof
By linearity it is sufficient to prove the associativity for the generators of $\lms$ in $T$:
\be\label{tras}
tr(x_\as, x_\bb \jp x_\gh) = tr(x_\as\jp x_\bb , x_\gh) \ ,\ \as,\bb,\gh \in T
\ee

We keep using the notation $\rho_1:= k_1+k_N$, $\rho_2:= k_1-k_N$, $\rho_3:= -k_2-k_3$.\\

Notice that $x_\as\jp( x_\bb \jp x_\gh)=0$ does not imply $(x_\as\jp x_\bb) \jp x_\gh=0$.\\
Let $\as = k_1+ k_4$, $\bb= k_1-k_4$, $\gamma = k_1-k_4$ then $x_\as\jp( x_\bb \jp x_\gh)=0$ (since  $x_\bb \jp x_\gh=0$) whereas $(x_\as\jp x_\bb) \jp x_\gh\ne 0$. However $tr(x_\as, x_\bb \jp x_\gh) = tr(x_\as\jp x_\bb , x_\gh)=tr(x_\gh)=0$.\\

We prove the associativity of the trace by checking case by case when $\as$ (or $\bb$ or $\gh$) is a {\it scalar P} ($\xu$ or $\xd$ or $\xt$) or a vector $V$ or a spinor $S$. We say $P\jp V\in V$ to mean that the $\jp$ product of a scalar and a vector is a vector, and similarly for other combinations of $P,V,S$. We also respectively denote by $TR1$ and $TR2$ the left-hand and right-hand sides of \eqref{tras}.

\bit
\item[PPP)] $tr(x_{P_i}, x_{P_j} \jp x_{P_\ell})= \delta_{j\ell}tr(x_{P_i} \jp x_{P_j}) = \delta_{ij}\delta_{j\ell} = tr(x_{P_i}\jp x_{P_j} , x_{P_\ell})$.
\item[PPV)] Since $(P\jp P)\jp V$ and $P\jp (P\jp V)$ are either $0$ or vectors then $TR1=TR2=0$
\item[PPS)] Since $(P\jp P)\jp S$ and $P\jp (P\jp S)$ are either $0$ or spinors then $TR1=TR2=0$
\item[PVV)]  Because of Lemma \sref{xpv},$TR1=TR2=0$ unless one vector is $x^\pm_{v\mu}$ and the other is  $x^\mp_{v\mu}$. We have $tr(x_{P_i}, x^\pm_{v\mu} \jp x^\mp_{v\mu})= \frac12 tr(x_{P_i} \jp (\xu+\xd)) = \frac12(\delta_{i1}+\delta_{i2})$; on the other hand $tr(x_{P_i}\jp x^\pm_{v\mu} , x^\mp_{v\mu})= \frac12 (\delta_{i1}+\delta_{i2})tr(x^\pm_{v\mu},x^\mp_{v\mu})=\frac12 (\delta_{i1}+\delta_{i2})\frac12 tr(\xu+\xd)=\frac12 (\delta_{i1}+\delta_{i2})$.
\item[PVS)] Since $(P\jp V)\jp S$, $P\jp (V\jp S)$, $V\jp (P\jp S)$ are either $0$ or spinors then $TR1=TR2=0$
\item[PSS)] We must show that
\be
tr(\xpi, x_{s_1} \jp x_{s_2}) = tr(\xpi\jp x_{s_1} , x_{s_2})
\ee
Since $x_{P_i}\jp x_{s_1} \propto x_{s_1}$, if $x_{s_1}\jp x_{s_2}=0$ then $TR1=TR2=0$.\\
Suppose now $x_{s_1}\jp x_{s_2}\ne 0$. We denote by $\rho_s$ the root associated to the spinor $x_s$. Then either $\rho_{s_1} + \rho_{s_2} = k_1-k_2-k_3 \pm k_N$ and
\be
 \rho_{s_1} -k_1\mp k_N + \rho_{s_2}  = -k_2-k_3
\ee
or, $\rho_{s_1} + \rho_{s_2} = k_1-k_2-k_3 \pm k_i$, for $4\le i \le N$ and
\be \rho_{s_1} + k_2 + k_3 + \rho_{s_2}  = k_1 \pm k_i
\ee
If $x_{s_1}\jp x_{s_2} = \lambda x_{vi}$ then $TR1=TR2=0$ and we are left with the case $\rho_{s_1} + \rho_{s_2} = k_1-k_2-k_3 \pm k_N$ hence either
$x_{s_1}\jp x_{s_2} = \frac12(\lambda_1 \xu +\lambda_3 \xt)$, which occurs if the coefficient of $k_N$ in both $\rho_{s_1}$ and $\rho_{s_2}$ is $(\rho_{s_1}, k_N)=(\rho_{s_2}, k_N)=+\frac12$, or $x_{s_1}\jp x_{s_2} = \frac12(\lambda_2 \xd +\lambda_3 \xt)$, which occurs if $(\rho_{s_1}, k_N)=(\rho_{s_2}, k_N)=-\frac12$.\\
Suppose $x_{s_1}\jp x_{s_2} = \frac12(\lambda_1 \xu +\lambda_3 \xt)$, then
\be
tr(\xpi, x_{s_1} \jp x_{s_2})= \frac12 tr(\xpi, \lambda_1 \xu +\lambda_3 \xt)= \frac12(\lambda_1 \delta_{i1} + \lambda_3 \delta_{i3})
\ee
on the other hand
\be
tr(\xpi\jp x_{s_1} , x_{s_2})= \frac12(\rho_i, \rho_{s_1})tr(x_{s_1}\jp x_{s_2}) = \frac14(\rho_i, \rho_{s_1})(\lambda_1+\lambda_3)
\ee
From table \ref{t:magiceo}, for $(r,s)=(1,1)$, if $\rho_i= \rho_1$ then $(\rho_i, \rho_{s_1})=1$; if $\rho_i= \rho_2$ then $(\rho_i, \rho_{s_1})=0$ and if $\rho_i= \rho_3$ then $(\rho_i, \rho_{s_1})=1$.\\
We calculate $\lambda_3$:
\bea{l}
-\ep(\rho_{s_1},-k_1-k_N)\ep(\rho_{s_1}-k_1-k_N, -\rho_{s_1} + k_1-k_2-k_3 + k_N)=\\
\ep(\rho_{s_1},-k_1-k_N)\ep(\rho_{s_1}-k_1-k_N, -k_2-k_3)=\\
\ep(\rho_{s_1},-k_1-k_N)\ep(\rho_{s_1}, -k_2-k_3)\ep(-k_1-k_N, k_2+k_3)
\eea
We calculate $\lambda_1$:
\bea{l}
-\ep(\rho_{s_1},k_2+k_3)\ep(\rho_{s_1}+k_2+k_3, -\rho_{s_1} + k_1-k_2-k_3 + k_N)=\\
\ep(\rho_{s_1},k_2+k_3)\ep(\rho_{s_1}+k_2+k_3, -k_1-k_N)=\\
\ep(\rho_{s_1},-k_1-k_N)\ep(\rho_{s_1}, -k_2-k_3)\ep(k_2+k_3,-k_1-k_N)
\eea
From $iv)$ of Proposition \sref{epsprop} together with \eqref{posroots1} and \eqref{posroots3}, we deduce that \be\label{uno}
\ep(-k_1\pm k_N, k_2+k_3)\ep(k_2+k_3,-k_1\pm k_N)=1
\ee
therefore $\lambda_1=\lambda_3$ from which $TR1=TR2$.\\
Similarly if $x_{s_1}\jp x_{s_2} = \frac12(\lambda_2 \xd +\lambda_3 \xt)$.
\item[VVV)] Since $(V\jp V)\jp V$ is either $0$ or a vector then $TR1=TR2=0$
\item[VVS)] Since $(V\jp V)\jp S$ and $V\jp (V\jp S)$ are either $0$ or spinors then $TR1=TR2=0$
\item[VSS)] This case is proven in the next Lemma \sref{VSS}.
\item[SSS)] Since $(S\jp S)\jp S$ is either $0$ or a spinor then $TR1=TR2=0$
\eit
\esquare

\begin{lemma} \label{VSS}
$tr(x_{vi}\jp x_{s1},x_{s2})= tr(x_{vi}, x_{s1} \jp x_{s2})$
\end{lemma}
\proof Let $\rho_{vi} = k_1\pm k_i$ and suppose $tr(x_{vi}\jp x_{s1},x_{s2})\ne 0$, then necessarily $\rho_{s1} = \frac12(k_1-k_2-k_3\pm ...\mp k_i \pm ... \pm k_N)$, namely the coefficients of $k_i$ must have opposite signs.
From table \ref{t:magiceo}:
\be\label{roo}
x_{so}:=x_{vi}\jp x_{s1} \imp \rho_{so} = \frac12(k_1-k_2-k_3\pm ...\pm k_i \pm ... \mp k_N) = \rho_{s1} \pm k_i \mp k_N
\ee
In order to have $tr(x_{vi}\jp x_{s1},x_{s2})\ne 0$ then, necessarily, either
\be
\rho_{so}-k_1\pm k_N + \rho_{s2} =-k_2-k_3
\ee
or
\be
\rho_{so}+k_2+k_3 + \rho_{s2} =k_1\mp k_N
\ee
In either case
\be\label{susd}
\rho_{so} + \rho_{s2} =k_1\mp k_N -k_2-k_3 \imp \rho_{s1} + \rho_{s2} =k_1 -k_2-k_3 \mp k_i
\ee
and
\bea{l}
tr(x_{vi}\jp x_{s1},x_{s2}) = \frac14 \ep(k_1\pm k_i, -k_1\mp k_N)\ep(\pm k_i \mp k_N, \rho_{s1})\times \\
\left[ \ep(\rho_{s1} \pm k_i \mp k_N, -k_1\pm k_N)\ep(\rho_{s1} -k_1 \pm k_i, -\rho_{s1}+ k_1 -k_2-k_3 \mp k_i) + \right. \\
\left. \ep(\rho_{s1} \pm k_i \mp k_N, k_2 + k_3)\ep(\rho_{s1} +k_2+k_3 \pm k_i \mp k_N, -\rho_{s1}+ k_1 -k_2-k_3 \mp k_i) \right]
\eea
Since $\rho_{s1} -k_1 \pm k_i$ and $\rho_{s1} +k_2+k_3 \pm k_i \mp k_N$ are roots then, using $i)$ of Proposition \sref{remas} and \eqref{uno}:
\bea{l}
tr(x_{vi}\jp x_{s1},x_{s2}) = -\frac14 \ep(k_1\pm k_i, -k_1\mp k_N)\ep(\pm k_i \mp k_N, \rho_{s1})\times \\
\left[ \ep(\rho_{s1} \pm k_i \mp k_N, -k_1\pm k_N)\ep(\rho_{s1} -k_1 \pm k_i, -k_2-k_3) + \right. \\
\left. \ep(\rho_{s1} \pm k_i \mp k_N, k_2 + k_3)\ep(\rho_{s1} +k_2+k_3 \pm k_i \mp k_N, k_1 \mp k_N) \right] =\\
-\frac14 \ep(k_1\pm k_i, -k_1\mp k_N)\ep(\pm k_i \mp k_N, \rho_{s1})\times \\
\ep(\rho_{s1} \pm k_i \mp k_N, -k_1\pm k_N)\ep(\rho_{s1}, -k_2-k_3)\times \\
\left[\ep(k_1 \mp k_i, k_2+k_3) + \ep(\pm k_i \mp k_N,k_2+k_3)\ep(k_2+k_3, k_1 \mp k_N)\right]=\\
-\frac14 \ep(k_1\pm k_i, -k_1\mp k_N)\ep(\pm k_i \mp k_N, \rho_{s1})\times \\
\ep(\rho_{s1} \pm k_i \mp k_N, -k_1\pm k_N)\ep(\rho_{s1}, -k_2-k_3)\times \\
\left[\ep(k_1 \mp k_i, k_2+k_3) + \ep(k_1 \mp k_i,k_2+k_3)\right]=\\
-\frac12 \ep(k_1\pm k_i, -k_1\mp k_N)\ep(\pm k_i \mp k_N, \rho_{s1})\times \\
\ep(\rho_{s1} \pm k_i \mp k_N, -k_1\pm k_N)\ep(\rho_{s1}, -k_2-k_3)\ep(k_1 \mp k_i, k_2+k_3)
\eea

Let us denote $x_{s1}\jp x_{s2}=\lambda_+ x_{s+}+\lambda_- x_{s-}+\lambda x_{v}$. From \eqref{susd}:
\bea{l}
\rho_{s\pm} = \rho_{s1} + \rho_{s2} -k_1 \pm k_N =-k_2-k_3 \mp k_i\pm k_N\\
\rho_{v} = \rho_{s1} + \rho_{s2} +k_2 + k_3 =k_1 \mp k_i
\eea
Since $\rho_{s\pm}$ is not a root then $x_{s\pm}=0$ and $x_{s1}\jp x_{s2} = \lambda x_v$. Together with $iii)$ of Lemma \sref{xpv} this implies
\bea{l}
tr(x_{vi},x_{s1}\jp x_{s2})=\lambda =\\ -\um \ep(\rho_{s1}, k_2 + k_3)
\ep (\rho_{s1} + k_2 + k_3, -\rho_{s1}+ k_1 -k_2-k_3 \mp k_i)(-1)^{\delta_{i,N-1}}=\\
\um \ep(\rho_{s1}, k_2 + k_3)
\ep (\rho_{s1} + k_2 + k_3, k_1 \mp k_i)(-1)^{\delta_{i,N-1}}
\eea

Therefore, in order to prove associativity, we must prove that:
\bea{l}\label{eq1}
-\ep(k_1\pm k_i, -k_1\mp k_N)\ep(\pm k_i \mp k_N, \rho_{s1})\times \\
\ep(\rho_{s1} \pm k_i \mp k_N, -k_1\pm k_N)\ep(k_1 \mp k_i, k_2+k_3) =\\
\ep (\rho_{s1} + k_2 + k_3, k_1 \mp k_i)(-1)^{\delta_{i,N-1}}
\eea
We have $\ep(\pm k_i \mp k_N, \rho_{s1}) = - \ep(\rho_{s1},\pm k_i \mp k_N)$ since, from \eqref{roo}, $\rho_{s1}\pm k_i \mp k_N$ is a root. Therefore we can write the left-hand side as:
\be
\ep(k_1\pm k_i, -k_1\mp k_N)\ep(\rho_{s1}, -k_1\pm k_i)\ep(\pm k_i \mp k_N, -k_1\pm k_N)\ep(k_1 \mp k_i, k_2+k_3)
\ee
Moreover $\ep(k_1 \mp k_i, k_2+k_3)= \ep(k_2+k_3,k_1 \mp k_i)$, by  $iv)$ of Proposition \sref{epsprop}, and the equality \eqref{eq1} reduces to
\be\label{dinu}
\ep(k_1\pm k_i, -k_1\mp k_N)\ep(\pm k_i \mp k_N, -k_1\pm k_N) =(-1)^{\delta_{i,N-1}}\\
\ee
By $iii)$ of Proposition \sref{epsprop}, we have $\ep(k_1\pm k_i, -k_1\mp k_N)=-\ep(-k_1\mp k_N ,k_1\pm k_i)$.
By adding and subtracting $-k_1\mp k_N$ on the right-hand side and, using $i)$ of Proposition \sref{remas} together with $vii)$ of Proposition \sref{epsprop}, we get
$\ep(k_1\pm k_i, -k_1\mp k_N)= \ep(-k_1\mp k_N, \pm k_i\mp k_N)$. The left-hand side of \eqref{dinu} becomes:
\bea{l}
\ep(k_1\pm k_i, -k_1\mp k_N)\ep(\pm k_i \mp k_N, -k_1\pm k_N)=\\
\ep(-k_1\mp k_N, \pm k_i\mp k_N)\ep(\pm k_i \mp k_N, -k_1\pm k_N)=\\
-\ep(-k_1\mp k_N, \pm k_i\mp k_N)\ep(-k_1\pm k_N,\pm k_i \mp k_N)=
-\ep(2k_1,\pm k_i \mp k_N)=\\
-\ep(\as_{N-1}-\as_{N-2}, \pm k_i \mp k_N) = (-1)^{\delta_{i,N-1}}
\eea
The last two equalities can be easily proven using \eqref{kal} and \eqref{posroots3}.\\
We have therefore proven that $tr(x_{vi}\jp x_{s1},x_{s2})= tr(x_{vi}, x_{s1} \jp x_{s2})$ if $tr(x_{vi}\jp x_{s1},x_{s2})\ne 0$.\\

Finally suppose $tr(x_{vi}\jp x_{s1},x_{s2})= 0$.\\
If $\rho_{s1} = \frac12(k_1-k_2-k_3\pm ...\pm k_i \pm ... \pm k_N)$ then $tr(x_{vi}\jp x_{s1},x_{s2})= 0$, as already shown.  In this case the vector $x_{s1}\jp x_{s2}$ has component $k_i$ either $0$ or $\pm 1$ and, from $iv),v)$ of Lemma \sref{xpv} also $tr(x_{vi},x_{s1}\jp x_{s2})=0$.\\
If $\rho_{s1} = \frac12(k_1-k_2-k_3\pm ...\mp k_i \pm ... \pm k_N)$ then we must have  $\rho_{s1} + \rho_{s2} \ne k_1 -k_2-k_3 \mp k_i$, otherwise $tr(x_{vi}\jp x_{s1},x_{s2})\ne 0$. Then neither $\rho_{s1} + \rho_{s2} -k_1\pm k_N= k_1\mp k_i$ nor $\rho_{s1} + \rho_{s2} +k_2 + k_3= k_1\mp k_i$
hence, by Lemma \sref{xpv}, $tr(x_{vi},x_{s1}\jp x_{s2})=0$. \esquare


\end{document}